\documentclass{amsart}
\usepackage{epsfig}

\begin{document}

\newtheorem{thm}{Theorem}[section]
\newtheorem{lem}[thm]{Lemma}
\newtheorem{cor}[thm]{Corollary}
\newtheorem{conj}[thm]{Conjecture}

\theoremstyle{definition}
\newtheorem{defn}{Definition}[section]

\theoremstyle{remark}
\newtheorem{rmk}{Remark}[section]
\newtheorem{exa}{Example}[section]

\def\square{\hfill${\vcenter{\vbox{\hrule height.4pt \hbox{\vrule width.4pt
height7pt \kern7pt \vrule width.4pt} \hrule height.4pt}}}$}

\def\R{\mathbb R}
\def\Z{\mathbb Z}
\def\CP{\mathbb {CP}}
\def\H{\mathbb H}
\def\E{\mathbb E}
\def\F{\mathcal F}
\def\D{\mathcal D}
\def\C{\mathbb C}
\def\til{\widetilde}
\def\N{\mathbb N}
\def\T{\mathcal T}
\def\G{\mathcal G}
\def\a{{\text{amb}}}

\newenvironment{pf}{{\it Proof:}\quad}{\square \vskip 12pt}

\title{Distortion of leaves in product foliations}

\author{Danny Calegari}
\address{Department of Mathematics \\ UC Berkeley \\ Berkeley, CA 94720}
\email{dannyc@math.berkeley.edu}

\maketitle

\begin{abstract}
We produce examples of codimension one foliations of $\H^2$ and $\E^2$ with
bounded geometry which are topologically product foliations, but for
which leaves are {\em non--recursively distorted}. That is, the function
which compares intrinsic distances in leaves with extrinsic distances in
the ambient space grows faster than any recursive function.
\end{abstract}

\section{Introduction}
It is a basic problem, given a foliation of a Riemannian manifold, to
compare the extrinsic and the intrinsic geometry of leaves of the foliation.
For leaves which are compact, this amounts to studying the distortion of
$\pi_1$ of the leaf in $\pi_1$ of the ambient manifold; much is known
about this problem. For a partial survey, see \cite{mG93}.

For taut codimension one foliations of $3$--manifolds, one knows by a
result of Sullivan in \cite{dS79} that one can choose a metric on the
ambient manifold such that the leaves are minimal surfaces. One could think
of this as saying that taut foliations ``measure area well''. However,
leaves of taut foliations are far from being quasi--isometrically
embedded, in general. For instance, a theorem of Fenley in \cite{sF92}
states that for an $\R$--covered foliation $\F$ of a hyperbolic $3$--manifold 
$M$, leaves of the pulled--back foliation $\til{\F}$ of $\til{M}$ limit to the
entire sphere at infinity of $\til{M} = \H^3$. A foliation is $\R$--covered
if the pulled--back foliation $\til{\F}$ is topologically the product foliation
of $\R^3$ by horizontal $\R^2$'s. Fenley's proof is attractive but 
somewhat lengthy, and therefore we think it is useful to give a 
simple proof of this fact here:

\begin{thm}
For $\F$ an $\R$--covered foliation of $M$, a finite volume 
hyperbolic $3$--manifold, every
leaf $\lambda$ of $\til{\F}$ limits to the entire sphere at infinity
$S^2_\infty$ of $\til{M} = \H^3$.
\end{thm}
\begin{pf}
For a leaf $\lambda$ of $\til{\F}$, let
$\lambda_\infty = \bar{\lambda} - \lambda \subset S^2_\infty$ be the
set of points in $S^2_\infty$ in the closure of $\lambda$. We claim that
$\lambda_\infty = S^2_\infty$ for each $\lambda$.

Suppose some $\lambda_\infty$ omits some point $p \in S^2_\infty$, and
therefore an open disk $U$ containing $p$. Let $U'$ and $U''$ be two
open disks in $S^2_\infty$ whose closures are disjoint and contained in $U$. 
Since $M$ has finite volume, there are elements $\alpha,\beta \in \pi_1(M)$
such that $\alpha(S^2_\infty - U) \subset U'$ and 
$\beta(S^2_\infty - U) \subset U''$. Let $\lambda' = \alpha(\lambda)$ and
$\lambda'' = \beta(\lambda)$. Since the leaf space of $\til{\F}$ is $\R$,
After relabeling the leaves if necessary, 
we can assume $\lambda$ separates $\lambda'$
from $\lambda''$. However, there is clearly an arc $\gamma$ in $S^2_\infty$
joining $\lambda'_\infty$ to $\lambda''_\infty$ which avoids $\lambda_\infty$.
$\gamma$ is a Hausdorff limit (in $\H^3 \cup S^2_\infty$) of arcs
$\gamma_i \subset \H^3$ running between $\lambda'$ and $\lambda''$. Each
$\gamma_i$ must intersect $\lambda$ in some point $q_i$, and by extracting
a subsequence, we find $q_i \to q \in \gamma$. But by construction,
$q \in \lambda_\infty$, giving us a contradiction.
\end{pf}

\vfill

If we want to study the distortion of leaves in taut foliations, then
$\R$--covered foliations are the most delicate case. For, if $\F$ is a
taut foliation which is {\em not} $\R$--covered, then the leaf space of
$\til{\F}$ is non--Hausdorff. That is, there are a sequence of pairs of
points $p_i,q_i \in \til{M}$ with $p_i \to p$ and $q_i \to q$ such that
$p_i,q_i \in \lambda_i$ but $p$ and $q$ lie on {\em different} leaves of
$\til{\F}$. It follows that the distance between $p_i$ and $q_i$ {\em as
measured in $\F$} is going to infinity; that is, 
$d_{\lambda_i}(p_i,q_i) \to \infty$. On the other hand, the distance between
$p_i$ and $q_i$ {\em as measured in $\til{M}$} is bounded:
$d_{\til{M}}(p_i,q_i) \to d_{\til{M}}(p,q)$.

By contrast, leaves of $\R$--covered foliations are never infinitely
distorted. In fact in \cite{dC99} we show the following:
\begin{thm}
Let $\F$ be an $\R$ covered foliation of an atoroidal $3$--manifold $M^3$.
Then leaves $\lambda$ of $\til{\F}$ are quasi--isometrically
embedded in their $t$--neighborhoods, for any $t$. That is, there is a
uniform $K_t,\epsilon_t$ such that for any leaf $\lambda$ of $\til{\F}$,
the embedding $\lambda \to N_t(\lambda)$ is a $(K_t,\epsilon_t)$ 
quasi--isometry.
\end{thm}

A heuristic scheme to measure the distortion of a leaf $\lambda$ runs
as follows.

We assume that $\F$ is co--oriented, so that in the universal cover there
is a well--defined notion of the space ``above'' and the space ``below''
a given leaf. Since leaves are minimal surfaces, the mean curvature is
zero, so there will be well--defined approximate directions at every
point in which the leaf ``bends upwards'' and approximate directions in
which it ``bends downwards''. For an arc $\alpha$ between points $p,q$ of
a leaf $\lambda$ which is roughly parallel to a direction of positive 
extrinsic curvature, the geodesic $\alpha^+$ 
in $\til{M}$ running between $p$ and $q$ will lie mostly above $\lambda$.
We can cap off the circle $\alpha \cup \alpha^+$ with a disk $D_\alpha$ of
minimal area, and consider its projection to $M$. The intersection of 
$D_\alpha$ with $\til{\F}$ gives a foliation of $D_\alpha$ by arcs parallel
to the principal directions of positive curvature in the leaves of $\til{\F}$.
Thus the projection of $D_\alpha$ to $M$ should be close to an embedding,
since its self--intersections cannot have a large dihedral angle.
If we consider longer and longer leafwise geodesics $\alpha$, the geometric
limit of the disks $D_\alpha$ in $M$ should be a geodesic lamination 
$\Lambda^+$ in $M$
which describes the ``eigendirections'' of positive curvature of $\F$.
If we consider the principal directions of negative curvature, we should get
a complementary lamination $\Lambda^-$, transverse to $\F$ and to $\Lambda^+$,
which describes the ``eigendirections'' of negative curvature of $\F$.

To study the distortion of leaves of $\F$ in such a setup, one need only 
study, for a typical leaf $l$ of $\Lambda^+$ say, how distorted the 
induced one--dimensional
foliation $\F \cap l$ is in $l$. Such a foliation has {\em bounded
geometry} --- its extrinsic curvature is bounded above and below by some
constant --- because of the compactness of $M$. Moreover it is topologically
a product. The point of this paper is to show that these two conditions in
no way allow one to establish any bound on the distortion function of
$\F$.

One remark worth making is that a pair of laminations $\Lambda^\pm$ 
transverse to an $\R$--covered foliation of an atoroidal $3$--manifold
$M$ are constructed in \cite{dC99} and independently by S\' ergio Fenley
(\cite{sF99}). The interpretation of these foliations as eigendirections
of extrinsic curvature is still conjectural, however.

\section{Non--recursive distortion}
\begin{defn}
The {\em distortion function} $\D$ of a foliation $\F$ in a Riemannian manifold
$M$ is defined as follows: if $d_\F$ denotes the intrinsic distance function
in a leaf of $\F$ thought of as a geodesic metric space, 
and $d_\a$ denotes the Riemannian distance function
in the ambient space, then
$$\D(t) = \sup_{x,y \in \lambda| d_\a(x,y) = t} d_\F(x,y)$$
where the supremum is taken over all pairs of points in all leaves of $\F$.
We usually expect that $M$ is simply connected.
\end{defn}

\begin{thm}
Let $\lambda$ be a smooth oriented bi--infinite ray properly 
immersed in $\H^2$.
\begin{enumerate}
\item{If $-1 \le \kappa \le 1$ everywhere, 
then the distortion function $\D$ is at most exponential.}
\item{If $\kappa \ge 1$ everywhere and $>1$ somewhere, 
then $\lambda$ has a self--intersection.}
\item{If $\kappa =1$ everywhere then $\lambda$ is a horocircle and therefore
has exponential distortion.}
\end{enumerate}
\end{thm}
\begin{pf}
The proof of 1. is a standard comparison argument, and amounts to 
showing that $\lambda$ makes at least as much progress as a horocircle.

To see that 2. holds, we consider the progress of the {\em osculating
horocircle} to $\lambda$. At a point $x$ on $\lambda$, the osculating
horocircle $H_x$ is the unique horocircle on the positive side of $\lambda$
which is tangent to $\lambda$ at $x$. Let $B_x \in S^1_\infty$ be the 
basepoint of $H_x$. Then the curvature condition on $\lambda$ implies that as
one moves along $\lambda$ in the positive direction, $B_x$ always moves 
anticlockwise. Since $\kappa>1$ somewhere, the derivative of $B_x$ is
nonzero somewhere. If we truncate $\lambda$ on some sufficiently big 
compact piece and fill
in the remaining segments with horocircles, 
we get a properly immersed arc $\lambda'$ in
$B^2$ with positive winding number, which consequently has a 
self--intersection somewhere. If this intersection is in the piece
agreeing with $\lambda$, we are done. Otherwise, the curvature condition
implies that $\lambda$ must have an intersection in the region enclosed
by $\lambda'$.
\end{pf}

In light of this theorem, it is perhaps surprising that we can make the
following construction:

\begin{thm}
For any $\epsilon>0$ there is a foliation $\F$ of $\H^2$ which is topologically
the standard foliation of $\R^2$ by horizontal lines with the following
properties:
\begin{itemize}
\item{Leaves are smooth.}
\item{The extrinsic curvature of any leaf at any point is bounded between
$1 - \epsilon$ and $1 + \epsilon$.}
\item{The distortion function $\D$ grows faster than any recursive function.}
\end{itemize}
\end{thm}
\begin{pf}
We consider the upper half-space model of $\H^2$ and let $\lambda$ be the
leaf passing through the point $i$. $\lambda$ will be the graph of a function
$r = \phi(\theta)$ in polar co--ordinates for $\theta \in (0,\pi)$. Choose
some very small $\delta$. Then for $\pi-\delta>\theta>\delta$ 
we let $\lambda$ agree with the horocircle with ``center'' at $\infty$ 
passing through $i$. Let $r_n$ be a sequence of positive real numbers which
grows faster than any recursive function. Then we define 
$$\phi \biggl( \frac \delta {2^n} \biggr) = r_n$$
Since $r_n$ grows so very fast, one can easily choose $\phi$ to interpolate
between $\frac \delta {2^n}$ and $\frac \delta {2^{n+1}}$ so that $\lambda$
is very close to a radial (Euclidean) line of very small slope in this range.
The extrinsic curvature of such a line in $\H^2$ is very close to $1$. 
Define $\phi$ in the range $\theta > \pi - \delta$ by 
$\phi(\theta) = \phi(\pi - \theta)$. Then the extrinsic distance between 
$(r,\theta)$ and $(r,\pi - \theta)$ depends only on $\theta$. It follows
that for $\theta$ between $\frac \delta {2^n}$ and $\frac \delta {2^{n+1}}$
the extrinsic distance between $(r,\theta)$ and $(r,\pi - \theta)$ is bounded
above by some recursive function in $n$. However, the length of the
path in $\lambda$ between $(r,\theta)$ and $(r,\pi - \theta)$ is growing
approximately like $\ln(r_n)$. The distortion function $\D$ is therefore
non--recursive.

Let $\lambda_1 = \lambda$ and for $t>0$, 
define $\lambda_t$ to be the graph of the
function $r = t \phi(\theta)$ in polar co-ordinates. Then there is a
hyperbolic isometry taking $\lambda_t$ to $\lambda_s$ for any $s,t$ and
therefore each of the $\lambda_t$ satisfies the same curvature bounds as
$\lambda$. Moreover, the union of the $\lambda_t$ gives a product foliation
of $\H^2$, as required.
\end{pf}

\begin{rmk}
For those unfamiliar with the concept, it should be noted that it is easy
to produce a function (even an integer valued function) such as $r_n$ which
grows faster than any recursive function. For, enumerate all the recursive
functions somehow as $\phi_n$. Then define $r(n) = \max_{m\le n} \phi_m(n)$.
Such an $r$ grows (eventually) at least as fast as any recursive function,
and therefore faster than any recursive function.
\end{rmk}

Essentially the same construction works for Euclidean space, a fact which we
now establish.

\begin{thm}
For any $\epsilon > 0$ there is a foliation $\F$ of $\E^2$ which is
topologically the standard foliation of $\R^2$ by horizontal lines with the
following properties:
\begin{itemize}
\item{Leaves are smooth.}
\item{The extrinsic curvature of any leaf at any point is bounded by 
$\epsilon$}
\item{The distortion function $\D$ grows faster than any recursive function.}
\end{itemize}
\end{thm}
\begin{pf}
The leaf $\lambda$ of $\F$ passing through the origin will be the graph of
an even function $y = \phi(x)$ and every other leaf will be a translate
of $\lambda$. After choosing a sufficiently large $K$ and small $\delta$,
we make $\phi(x) = \delta x^2$ for $|x| < K$, so that $\lambda$ has
curvature bounded below $\epsilon$, and is very nearly vertical at $x = K$.
Then let $r_i$ be a sequence growing faster than any recursive function,
as before, and define $\phi(K+n) = r_n$. Then we can choose $\phi$ to
interpolate between $K+n$ and $K+n+1$ for each $n$ to make it very smooth
and almost straight, since $\lambda$ will be almost vertical in this region.
\end{pf}

\section{Minimality}

\begin{defn}
A {\em decorated} metric space is a metric space in the usual sense
with some auxiliary structure. For instance, this structure could consist
of a basepoint, some collection of submanifolds, a foliation or lamination, 
etc. together with a notion of convergence of such structures 
in the geometric topology, in the sense of Gromov.
One says a decorated metric space $M$ has {\em bounded geometry} if the metric
space itself has bounded geometry, and if for every $t$, the decorated
metric spaces obtained as restrictions of $M$ to the balls of radius $t$
form a precompact family.
\end{defn}

For instance, we might have a family of foliated Riemannian manifolds
$M_i,\F_i$. The geometric topology requires a choice of basepoints $p_i$
in each $M_i$. Then we say the family $(M_i,\F_i,p_i)$ is a Cauchy sequence
if there are a sequence of radii $r_i$ and $\epsilon_i \to 0$ so that
the Gromov--Hausdorff distance between the ball of radius $r_i$ in
$M_i$ and $M_j$ about $p_i,p_j$ 
for $j>i$ is at most $\epsilon_i$, and that such
``near--isometries'' can be chosen in a way that the leaves of $\F_i$
can be taken $\epsilon_i$--close to the leaves of $\F_j$. From such
a Cauchy sequence we can extract a limit $(M,\F,p)$.

\begin{defn}
A decorated metric space $M$ with bounded geometry is 
{\em minimal} if for every limit $(M^+,p)$
of the pointed decorated metric spaces $(M,p_i)$ and every point $q \in M$, 
there is a sequence $q_i \in M^+$ such that the pointed decorated metric
spaces $(M^+,q_i) \to (M,q)$.
\end{defn}

The foliations constructed in the last section are certainly not minimal.
We ask the question: does the assumption of minimality allow one to 
make some estimate on the distortion function for a product foliation?

\end{document}